\documentclass{article}
\usepackage{amsthm,amsfonts,amsbsy, amssymb,amsmath,graphicx}
\usepackage{graphics}

\newtheorem{thm}{Theorem}
\newtheorem{lm}{Lemma}

\newtheorem{re}{Remark}
\newtheorem{prop}{Proposition}

\title{The Khovanov Complex for Virtual Links} \label{hvnv}
\date{}
\author{Vassily Olegovich Manturov}

\begin{document}

\maketitle

\section{Introduction}

\begin{figure}
\centering\includegraphics[width=360pt]{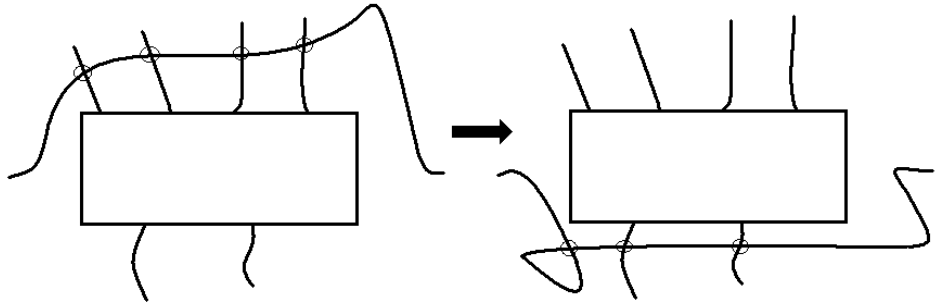} \caption{The
detour move} \label{obj}
\end{figure}

In the last few years, knot theory has enjoyed a rapidly
developing generalisation, the {\em Virtual knot theory}, proposed
by Louis Kauffman in 1996, see \cite{KaV}. A virtual link is a
combinatorial generalisation of the notion of classical links: we
consider {\em planar diagrams} with a new crossing type allowed;
this new crossing (called virtual and marked by a circle) is
neither an overcrossing nor an undercrossing. It should be treated
as an artefact of two branches, which do not want to intersect but
can not do without. This leads to a natural generalisation of
Reidemeister moves for the virtual case: besides usual ones (which
should be treated as local transformation inside a small
3-dimensional domain), we also add a {\em detour move}, which
means the following. If there is arc of the diagram between some
points $A$ and  $B$ contains only virtual crossings, it can be
removed and detoured as any other path connecting  $A$ and $B$;
all crossings which occur in this path are set to be virtual as
shown in Fig. \ref{obj}.

Thus, {\em a virtual link} is an equivalence class of virtual
diagrams modulo generalised Reidemeister moves, the latter
consisting of classical Reidemeister moves and the detour move.

\begin{re}
This definition allows to generalise some classical knot
invariance; herewith the invariance under the detour move often
happens to follow from the definition, whence the classical move
invariance repeats the invariance for the case of classical
knots.\label{rRe}
\end{re}

This approach leads to the following geometric definition for a
virtual link: it is a link in a thickened $2$-surface $S_{g}\times
I$ considered up to homotopy and stabilisation/destabilisation;
herewith virtual crossings correspond to ``projection artefacts''
while ``projecting'' the knot theory from $S_{g}\times I$ to ${\bf
R}^{2}\times I$ (here $S_{g}$ is projected to the plane, whence
the images of classical crossings represent classical crossings
and intersections of branches, which do not intersect in
$S_{g}\times I$, form virtual crossings. Thus, virtual links admit
a geometrical definition, besides the combinatorial one.

Many known classical knot invariants were generalised for the
virtual case, see, e.g. \cite{KaV,Ma,Ma6,Ma'9}, in particular the
Kauffman bracket and the Jones polynomial, \cite{KaV}.

An outstanding achievement of the recent years classical knot
theory of  is a generalisation of the Jones polynomial (Kauffman
bracket) proposed by Khovanov, \cite{Kh1}. This generalisation is
a {\em categorification} of the Jones polynomial in the following
sense. We replace abstract polynomials by (co)homologies of graded
complexes of vector spaces. Thus, with each link we accociate a
certain graded algebraic complex, such that all homologies of this
complex are knot invariants, and the graded Euler characteristic
coincides with the Jones polynomial.

\newcommand{\cC}{{\cal C}}

In the sequel, we shall consider graded complexes
$\cC=\oplus_{i,j}\cC^{i,j}$, with {\em height} $i$ and {\em
grading} $j$. The differential mapping is assumed to preserve the
grading and shift the height by $+1$.

It would be more reasonable to call such complexes {\em
cohomological}, but it is usual to say about {\em Khovanov
homologies}, thus we shall say {\em chains, cycles, boundaries}
instead of {\em cochains, cocycles, coboundaries}. For such
complexes we define the {\em height shift} $\cC\mapsto \cC[k]$ and
the {\em degree shift} $\cC\mapsto \cC\cdot\{l\}$ as follows:
$(\cC[k])^{i,j}=\cC^{i-k,j}$;
$({\cC\cdot\{l\}})^{i,j}={\cC}^{i,j-l}$. In the first case (height
shift) all differentials are shifted accordingly. By {\em graded
Euler characteristic} of a complex $\cC^{ij}$ is meant the
alternating sum of graded dimensions of chain spaces, or,
equivalently, that of homology spaces. For chain spaces we have
$\chi(C)=\sum_{i}(-1)^{i}qdim {\cC}^{i}=\sum_{i,j} (-1)^{i}q^{j}
dim {\cC}^{i,j}$.

An outstanding property of the Khovanov homology is {\em
functoriality} meaning the following. Having a cobordism in ${\bf
R}^{3}\times I$ between two links $K_{1}\subset {\bf
R}^{3}\times\{0\}$ и $K_{2}\subset {\bf R}^{3}\times \{1\}$, we
get a natural mapping between Khovanov homologies $Kh(K_{1})\to
Kh(K_{2})$, which represents a {\em cobordism invaraint} (up to
sign $\pm 1$); in particular, when links  $K_{1}$ и $K_{2}$ are
empty, the Khovanov homology gives an invariant of  {\em 2-links},
i.e. compact 2-surfaces embedded in  ${\bf R}^{4}$. This result
was obtained independently by Khovanov, Bar-Natan, and Jacobsson.
For more details see \cite{BN6}.

The aim of the present work is to generalise the Khovanov
construction for the case of virtual knot. In view of Remark
\ref{rRe}, the invariance proof for Khovanov homologies for
virtual links can be simplified. Here, the most difficult part is
the {\em definition} of the Khovanov complex to be really a graded
complex (with $d^{2}=0$). This work is represented in more details
in the book \cite{Ma}.

We shall describe three ways for constructing the Khovanov
homology for virtual links. The first one deals with arbitrary
links and Khovanov homologies with ${\bf Z}_{2}$-coefficients; the
second approach leads to the Khovanov homologies for framed links
(these homologies are constructed by doubling the initial link),
the third approach is based on two-fold orineting coverings over
virtual knots (in sense of atoms).

\section{The Kauffman Bracket and Jones Polynomial. \\ Atoms}

\begin{figure}
\centering\includegraphics[width=100pt]{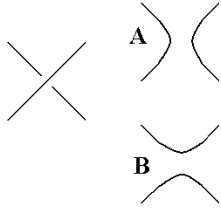} \caption{Two
ways of smoothing} \label{smth}
\end{figure}

Consider an oriented virtual diagram  $L$ and the diagram $|L|$
obtained from  $L$ by ``forgetting'' the orientation. Let us
smooth classical crossings of the diagram $|L|$ according to the
following rule. Each classical crossing can be {\em smoothed} in
one of two ways: $A$ и $B$, see Fig. \ref{smth}. The way of
smoothing for all classical crossings of a diagram generates a
{\em state} of it. Each state genrates a set of curves on the
plane; these curves have only virtual crossings. In other words,
we have got a trivial virtual links; its components are unicursal
curves.

Suppose the diagram $L$ has $n$ classical crossings. Let us
enumerate them arbitrarily. Thus we have got  $2^{n}$ states;
these states are in one-to-one correspondence of the $n$-cube
$\{0,1\}^{n}$, where $0$ and $1$ correspond to $A$-smoothings and
$B$-smoothing, respectively. Each state $s$ has three important
characteristics: the number  $\alpha(s)$ of  $A$-type smoothings,
the number $\beta(s)=n-\alpha(s)$ of $B$-type smoothings and the
number  $\gamma(s)$ of link components in the corresponding state.
After that, we define the Jones-Kauffman polynomial by the
following formula:

\begin{equation}X(L)=(-a)^{-3w(L)}\sum_{s}a^{\alpha(s)-\beta(s)}(-a^{2}-a^{-2})^{\gamma(s)-1},
\label{18n1}\end{equation} where $w(L)$ is the linking number of
the oriented diagram $L$ (i.e. the difference between the number
of positive classical crossings  and the number of negative
classical crossings), and $\alpha(s),\beta(s)$ are the numbers of
$A$-type and $B$-type smoothings, respectively.

The unnormalised version of the the Jones polynomial
$\sum_{s}a^{\alpha(s)-\beta(s))}(-a^{2}-a^{-2})$ is called {\em
the Kauffman bracket}. The Kauffman bracket is defined for
unoriented links; it is invariant under the detour move,
$\Omega_{2}, \Omega_{3}$ and the doubled version of the first
Reidemeister move $\Omega_{1}^{2}$ which consists in addition
(removal) of two ``adjacent'' loops, for more details, see
\cite{Ma}.

\newcommand{\Jj}{\hat J}

\begin{figure}
\centering\includegraphics[width=200pt]{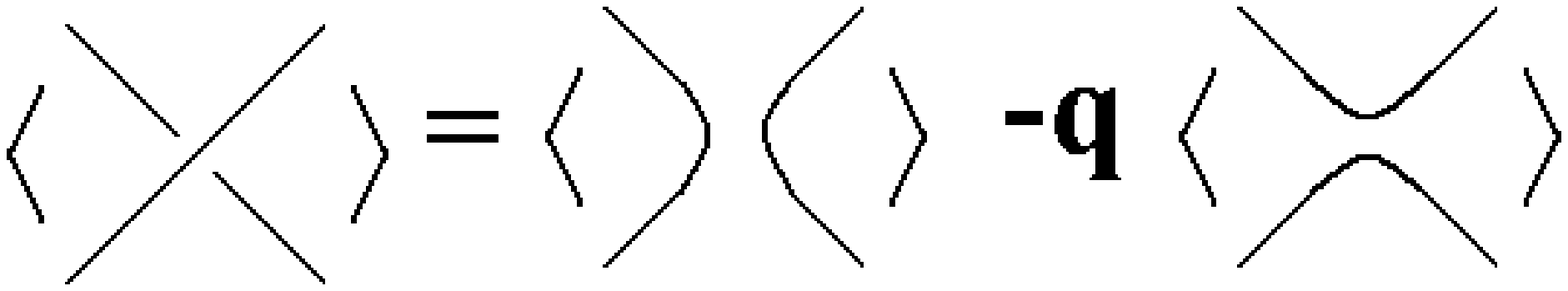} \caption{The
modified Kauffman relation} \label{smth2}
\end{figure}

The variable change $a=\sqrt{(-q^{-1})}$ transformed the Jones
polynomial to its modified version $J$. We shall also use a
version denoted by $\Jj$. They differ by normalisation; thus,
$J=1$ on the unknot, whence $\Jj=1$ on the empty unlink. Herewith,
$J=\frac{\Jj}{(q+q^{-1})}$. Let us describe this construction in
more details. Let $L$ be an oriented virtual diagram with $n$
classical crossings and let $|L|$ be the corresponding
non-orientable virtual diagram. Let $n_{+}$ and $n_{-}$ denote the
numbers of positive and negative crossings of $L$, so that
$n=n_{+}+n_{-}$. Put

\begin{equation}\Jj(L)=(-1)^{n_{-}}q^{n_{+}-2 n_{-}}\langle
L\rangle,\label{nenorm}\end{equation}

where $\langle L\rangle$ is the modified Kauffman bracket defined
according to the rule $\langle \bigcirc\rangle=1$,  and the
modified Kauffman relation, see Fig. \ref{smth2}.

In the sequel, we shall deal with $\Jj$.

The equivalence classes of virtual diagrams modulo the detour move
and the moves $\Omega_{1}^{2},\Omega_{2},\Omega_{3}$ are called
{\em framed virtual links}; for framed virtual links there is a
well defined operation of taking $n$ ``parallel copies'': $K\to
D_{n}(K)$: when applying one of the moves
$\Omega_{1}^{2},\Omega_{2},\Omega_{3}$ or the detour move  (but
not  $\Omega_{1}$) to the source diagram $K$, the diagram
$D_{n}(K)$ is transformed to an equivalent diagram .

The polynomial  $\Jj$ has a simple combinatorial description in
terms of the state cub. Taking off the normalising coefficients
$(-1)^{n_{-}}q^{n_{+}-2 n_{-}}$, we obtain a Kauffman bracket of
the type $\sum_{s}(-q)^{\beta(s)}((q+q^{-1})^{\gamma(s)})$, i.e.
we take the sum over all vertices of the cube of $(-q)$ to the
{\em height of the vertex} multiplied by $(q+q^{-1})$ to the
number of circles corresponding to the vertex. Thus we can say
that we replace each circle by $(q+q^{-1}$, then the polynomials
are multiplied and ``shifted'' by some $\pm q^{k}$.

This means that the Jones polynomial can be recovered only by the
information about the {\em numbers of vertices} in the cube state.
Taking into account the information, {\em how the circles are
reconstructed while passing from a state to an adjacent one}, we
obtain the Khovanov homology.

Now, we describe an important construction that we shall use in
the sequel. An {\em atom} is a pair $(M,\Gamma)$ : compact
$2$-manifold  $M$ without boundary and a $4$-valent graph
$\Gamma\subset M$ dividing $M$ into cells with a fixed
checkerboard colouring. The atoms are considered up to natural
combinatorial equivalence: homeomorphisms of the manifolds mapping
the frame to the frame and preserving the colour of cells. We
shall not assume $M$ to be connected; sometimes we shall need the
case when $M$ consists of two connected components.

Each atom (more precisely, its equivalence class) can be recovered
from the following combinatorial data:

\begin{enumerate}
\item The frame ($4$-valent graph);

\item The $A$--structure (specifying for each vertex, which
outgoing semiedges are {\em opposite}; this data is defined
according to the local structure of opposite edges on the surface)
and

\item $B$--structure (at each vertex we mark two pairs of adjacent
(not opposite) semiedges (or, equivalently, two pairs of opposite
angles) to form the local boundary of two black cells.
\end{enumerate}

Starting from 1996, the author has been developing the connection
between atoms and knots and virtual knots. We shall describe this
connection later in the text. For more details, see \cite{Ma}.

\section{The Khovanov Complex with ${\bf Z}_{2}$-Coefficients]}

In the present section, we present a result first proved in
\cite{Ma'9}.

Let $L$ be an orientable virtual link diagram. By {\em perestroika
cube} we mean the cube $\{0,1\}^{n}$ at each vertex of which we
indicate the number of state (as in the state cube), and for each
edge we show {\em how the set of circles is transformed} (new
information). Associate with each circle the linear space $V$
over~${\bf Z}_{2}$ generated by the two vectors $v_{+}$ and
$v_{-}$, such that the grading of $v_{\pm}$ is equal to $\pm 1$.
To each vertex $s=\{a_{1},\dots,a_{n}\}$ of the cube, there
corresponds some number of circles to be defined by $\gamma(s)$.
We replace such a vertex by the vector space $V^{\otimes
\gamma(s)}\cdot\{\sum_{i=0}^{n}a_{i}\}$ which is obtained from the
tensor power of the space  $V$ by a degree shift. This replacement
of $(q+q^{-1})$ by $V$  such that $\mbox{qdim}V=(q+q)^{-1}$ is a
very important step in the categorification. Now we define the
space of chains of the height $k$ as the tensor sum of spaces
related to all vertices of the height  $k$.

Now, let us define the partial differential on the chains which
act along the edges of the cube, as follows. Let an edge $a$
correspond to the switch from a state $s$ to a stat $s'$; herewith
the $l$ circles not adjacent to the crossing in question, are not
changed. At the crossing of the diagram $|L|$ corresponding to the
edge $a$, either one circle is transformed to two circles, or two
circles are transformed to one circle, or one circle is
transformed to one circle. In the first two cases we define the
differential as before, namely, $\Delta\otimes Id^{\otimes l}\cdot
\{1\}$ или $m\otimes Id^{\otimes l}\cdot\{1\}$. Here the mappings
$m:V\otimes V\to V$ and $\Delta:V\to V\otimes V$ are defined as
follows.

The mapping $m$:

\begin{equation}\left\{{v_{+}\otimes v_{-}\mapsto v_{-}, v_{+}\otimes v_{+}\mapsto v_{+},}\atop {v_{-}\otimes
v_{+}\mapsto v_{-}, v_{-}\otimes v_{-}\mapsto
0}\right.\end{equation}

The mapping $\Delta:$

\begin{equation}\left\{ { v_{+}\mapsto v_{+}\otimes v_{-}+v_{-}\otimes v_{+}} \atop{v_{-}\mapsto v_{-}\otimes
v_{-}.}
 \right.\end{equation}\label{mDe}

In the case of  $(1\to 1)$-perestroika, we define the ``partial
differential'' on the edge to map the whole space to zero. The
question how to define the differential of type $(1\to 1)$ well is
the main difficulty in the common case; in the case of the field
${\bf Z}_{2}$ this difficulty can beovercome. Denote the obtained
perstroika cube by $Q(L)$. In order for the differential to be
well defined, the cube should be {\em commutative}, that is, for
any 2-face of the cube, the composition of two mappings
corresponding to one sequence of edges, should be equal to the
(minus) composition corresponding to the other sequence of edges.
Note that in this case (for the field ${\bf Z}_{2}$) the
anticommutativity is equivalent to the commutativity.

\begin{lm}
The cube $Q(L)$ is commutative.
\end{lm}

\begin{figure}
\centering\includegraphics[width=150pt]{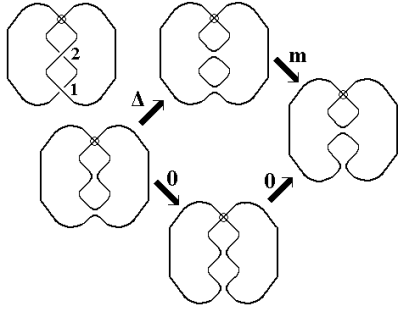}
\caption{Checking the commutativity of a $2$-face} \label{hovaz2}
\end{figure}

This statement follows from a routin check analogous to that
performed by Bar-Natan in \cite{BN4}. We shall check only one case
(the most interesting), see Fig. \ref{hovaz2}.

We have to show that the $2$-mapping $m\circ \Delta:V\to V$ maps
all space $V$ to zero. Indeed, we have $v_{-}\mapsto v_{-}\otimes
v_{-}\mapsto 0$, $v_{+}\mapsto v_{+}\otimes v_{-}+v_{-}\otimes
v_{+}\mapsto 2 v_{-}=0$ over the field~${\bf Z}_{2}$.

Note that this is the only essential non-classical case, where we
have a $1\to 1$-type perestroika. Indeed, it follows from pariry,
that the number of $1-1$-perestroikas does not exceed two. If
there are no such perestroikas at all, the question is reduced to
one of the classical cases (all considered by Bar-Natan).

If there were two or four such perestroikas, then the $2$-face in
question

\begin{equation}
\begin{array}{ccc} V^{\otimes a}\{1\} & \stackrel{s}{\longrightarrow} & V^{\otimes b}\{2\} \cr
r\uparrow & &\uparrow q \cr V^{\otimes c} &
\stackrel{p}{\longrightarrow} & V^{\otimes d}\{1\}
\end{array}
\end{equation}
either each of the two compositions $q\circ p$ and $s\circ r$
contains a zero mapping $1\to 1$ (for instance, in the case $a=b,
c=d$ the mappings $p$ and $s$ are zeroes), or we have the case
described above.

Put $\cC(L)=Q'(|L|)\{n_{+}-2n_{-}\}[n_{-}]$. In this case $\cC(L)$
is a well-defined complex. Denote the (bi)graded linear space of
homologies for $\cC(L)$ by $Kh(L)$ (or by $Kh_{{\bf Z}_2}(L)$ if
we wish to underline that the Khovanov homology is considered over
the field ${\bf Z}_{2}$).

\begin{thm}
$Kh(L)$ is an invariant of $L$. Furthermore, the graded Euler
characteristic $\chi(Kh(L))$ equals $\Jj(L)$.
\end{thm}

The invariance of the Khovanov homology under classical
Reidemeister moves is verbally the same as that proposed by
Bar-Natan \cite{BN4} for classical links; it is classical and
deals with a local domain where a Reidemeister move is performed
(and does not pay attention to what happens outside). The
invariance under the detour move follows from the definition.

The second statement of the theorem follows from the fact that the
Euler characteristic defined as the alternating sum of (graded)
Betti numbers is equal the alternating sum of graded dimensions of
chain spaces.

Note that the complex $\cC(L)$ is splitted into two complexes:
those of even grading and of odd grading (we remind that the
differential preserves the grading).

Thus, we get {\bf two} types of the Khovanov homology: the even
one $Kh^{e}$ and the odd one $Kh^{o}$.

They correspond to monomials of the Jones polynomial having
degrees congruent to two modulo four ($Kh^{e}$) and monomials with
degree divisible by four ($Kh^{o}$). For any classical link, we
have only one of these sorts of homologies. More precisely, the
following theorem holds.

\begin{thm}
For a classical even-component link, we have $Kh_{o}=0$; for a
classical odd-component link, we have $Kh_{e}=0$.
\end{thm}

\begin{figure}
\centering\includegraphics[width=180pt]{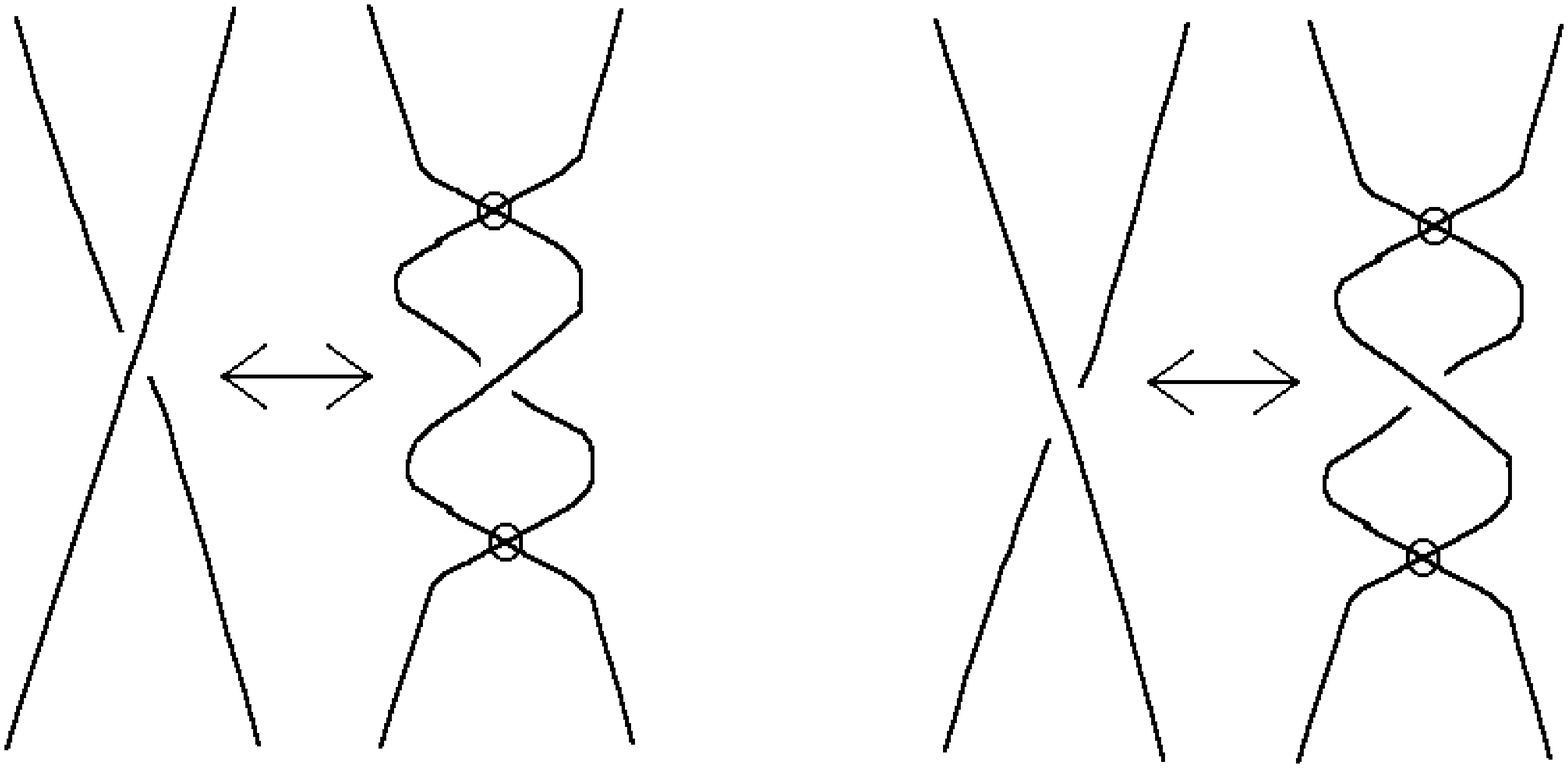} \caption{Two
versions of virtualisation} \label{twist}
\end{figure}

By a {\em virtualisation} we mean the following local
transformation in a neighbourhood of a classical crossing, see
Fig. \ref{twist}.

Note that the ${\bf Z}_{2}$-Khovanov complex, we have constructed
is completely defined by the perestroika cube and the numbers
$n_{+},n_{-}$. Thus, the Khovanov homologies are invariant under
virtualisations of the diagram.

In the next section, we shall give another construction of
Khovanov homologies (for framed links). This approach is
virtualisation-sensitive. The Khoavanov complex given in this
section coincides with the usual ${\bf Z}_{2}$-Khovanov complex in
the case of classical links; as mentioned above, it is not to
handle the $1\to 1$-differentials here. In the next two sections,
we shall consider the Khovanov complex not for all virtual links,
but only for ``good'' ones: those for which the cube has no local
differentiasl of the type $1\to 1$. Later on, we shall transform
each virtual link (diagram) to a good diagram and observe the
Khovanov homologies of the transformed diagram while performing
the Reidemeister moves to the initial diagram. Thus, in the the
section, we shall construct the Khovanov complex for framed links,
where the good diagram is the ``doubled'' version of the initial
diagram.

\section{The Khovanov Complex for Doubled Knots}

In the remaining part of the present chapter, we shall use the
following construction connecting atoms and virtual links.

Given an atom  $V$. Consider a generic embedding $p$ of its frame
in ${\bf R}^{2}$ preserving the $A$-structure of the atom. Thus we
get a $4$-valent graph on the plane with two types of vertices:
images of atom's vertices and ``immersion artifacts'', that is,
intersection of images of atom's edges. In the first case, we mark
crossings to be classical as in the case of classical knots
and height atoms.

Namely, let  $a,b$ be semiedges emanating from a vertex $X$ and
forming a black angle according to the $B$-structure of the atom.
While projecting the neighbourhood of $X$, denote by $u$ the one
of semiedges $a,b$, which is passed {\bf before} while sweeping
the angle $X,a,b$ in the clockwise direction. The semiedge  $u$
and the semiedge opposite to it are chosen to form an
undercrossing at $X$.

The intersection points of images of different edges are to be
marked by virtual crossings. Thus we obtain a virtual link diagram
to be denoted by $L$.

\begin{lm} The virtual link $L$ is defined by the atom up to
detour moves and virtualisations of some classical
vertices.\label{luma}
\end{lm}

\begin{proof} While projecting a graph to the plane we respect only the structure
of opposite semiedges. Suppose for some vertex $A$ of the atom in
question we have four emanating edges $a,b,c,d$, so that the edge
$c$ is opposite to $a$ and the edge $d$ is opposite to $b$. While
immersing this to ${\bf R}^{2}$ the cyclic ordering (in the
counterclockwise direction) may be either $a,b,c,d$ or $a,d,c,b$.
It is easy to see that the virtual links obtained in this way
differ from each other by a virtualisation at the corresponding
classical crossing (and detour moves).
\end{proof}

Let us mention one more importan.

\begin{lm}
Let  $L,L'$ be two virtual diagrams, such that the corresponding
atoms are orientable, and  $L'$ can be obtained from $L$ by using
a classical Reidemeister move. Then $Kh(L)=Kh(L')$.\label{lemm}
\end{lm}

While performing the detour move, the structure of classical
crossings is not change, thus we do not change the perestroika
cube and the whole complex. In the case of {\em classical
Reidemeister moves} there exists an invariance proof by Bar-Natan.
It is local, i.e., it uses only the local structure of the
Reidemeister move (regardless what happens outside). Thus, it
works for virtual knots as well.

\begin{prop}
Let $L$ be a virtual link diagram. Then the atom corresponding to
the diagram$ D_{2}(L)$ is orienable.\label{propp}
\end{prop}

For the proof, we shall need one more auxiliary proposition.

Given an atom $V$. Suppose its  $A$-structure is such that there
exists an orientation for all edges of  $V$ such that for each
vertex, some two opposite (semi)edges are outgoing, and the other
two are incoming. We call such an orienation of edges a {\em
source-target structure}.

\begin{prop}[\cite{Ma}]
An atom admitting a source-target structure is
orientable.\label{pris}
\end{prop}

\begin{figure}
\begin{center}
\unitlength 1.2pt
\begin{picture}(100,100)
\put(55,5){\vector(0,1){90}} \put(45,95){\vector(0,-1){90}}
\put(5,45){\line(1,0){38}} \put(47,45){\line(1,0){6}}
\put(57,45){\vector(1,0){38}} \put(95,55){\line(-1,0){38}}
\put(53,55){\line(-1,0){6}} \put(43,55){\vector(-1,0){38}}
\put(80,45){\vector(0,1){5}}
\put(80,51){\makebox(0,0)[cc]{\tiny{$\tau_{A}$}}}
\put(80,45){\circle*{1}}
\put(78,43){\makebox(0,0)[cc]{\tiny{$A$}}}

\put(87,47){\makebox(0,0)[cc]{\tiny{$\partial L_{A}$}}}
\end{picture}
\end{center}
\caption{Local orientation of the doubled diagram} \label{nizhe}
\end{figure}
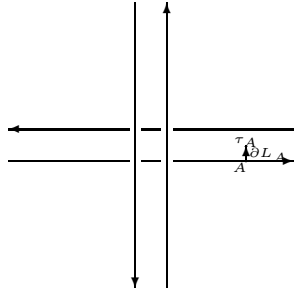

Then, the proof of Proposition \ref{propp} is deduced from
Proposition \ref{pris} as follows. Let $L$ be a virtual diagram.
Orient the diagram $D_{2}(L)$ in such a way that for any point $A$
on the diagram the reference frame $({\partial L}_{A},\tau_{A})$
generates a positive orientation of the plane, see Fig.
\ref{nizhe}.

Here $\partial L$ is the tangent vector for the link and
$\tau_{A}$ is the vector perpendicular to $\partial L$ directed
from the point $A$ to a close point on the adjacent component.

The desired source-target structure is now constructed as follows:
all edges (i.e., images of edges of the corresponding atom) of the
diagram $D_{2}(L)$ are naturally splitted into ``long'' ones
(those corresponding to edges of $L$) and ``short'' ones (four
shourt edges correspond to each vertex of $L$, see Fig.
\ref{nizhe}). Let us change the orientation of short edges without
changing that of the long ones. The obtained orienation gives a
source-target structure.

\begin{re}
This structure plays a key role in some other problems of
combinatorial topology, see, e.g. \cite{Izv}.
\end{re}

Thus, the Khovanov complex for $D_{2n}(L)=D_{2}(D_{n}(L))$ is well
defined for each ring of coefficients. The mapping  $L\mapsto
D_{2n}(L)$ is a framed link invariant. Thus, it is natural to
expect that the Khovanov homology of the doubled link are
invariants of the initial link. More precisely, the following
lemma takes place.

\begin{lm}
Let $L,L'$ be two equivalent framed virtual link diagrams. Then
there exists a sequence of diagrams $D_{2}(L)=L_{0},L_{1},\dots,$
$L_{n}=D_{2}(L')$ such that \label{lima}

\begin{enumerate}

\item All atoms corresponding to diagrams $L_{i}$ are orientable.

\item For each $i=0,\dots, n-1$, the diagram $L_{i+1}$ can be
obtained from  $L_{i}$ by a detour or by a Reidemeister move.
\end{enumerate}
\end{lm}

\begin{figure}
\centering\includegraphics[width=300pt]{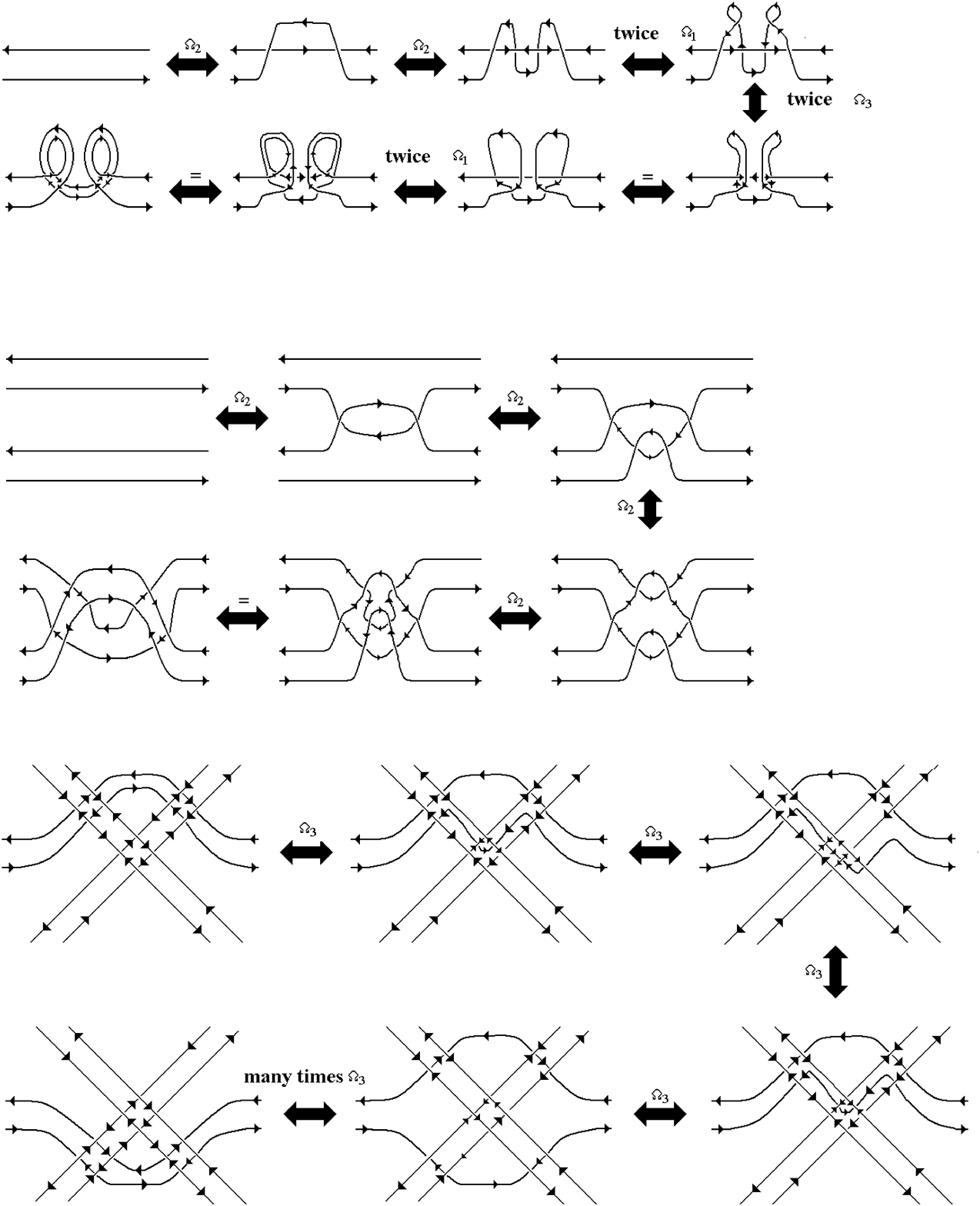}
\caption{Doubled versions of Reidemeister moves} \label{dvig''}
\end{figure}

The proof of Lemma \ref{lima} is shown in Fig. \ref{dvig''}.

\begin{thm}
Let $n$ be a positive integer. Then the image of the map $L\mapsto
Kh(D_{2n}(L))$ is a virtual link invariant.
\end{thm}

\begin{proof}
In view of Proposition \ref{propp}, $\cC(D_{2n}(L))$ is a well
defined complex. Now let $L,L'$ be two equivalent diagrams of
framed virtual links. Then, by lemma \ref{lima}, there exists a
sequence of diagrams $D_{2n}(L)=L_{0},\dots,L_{n}=D_{2n}(L')$ all
corresponding to orientable atoms such that $L_{i+1}$ is obtained
from $L_{i}$ by a generalised Reidemeister moves. Thus, in view of
Lemma \ref{lemm}, we have
$Kh(D_{2n}(L))=Kh(L_{1})=\dots=Kh(D_{2n}(L')$. The theorem is
proved.
\end{proof}

Note that the doubled diagram for $L$ and the doubled version for
the diagram $L'$ obtained from $L$ by virtualising a crossing
 of $L$ have quite different state cubes. Thus, the complex we
 have constructed in the present section, should be
 virtualisation-sensitive.

However the ``doubled'' Khovanov complex we have constructed
conceptually differs from the usual one known for the case of
classical links: here both in classical and virtual case, one
should first double the diagram, and then calculated the Khovanov
homologies.

A very natural question is whether the ``usual'' Khovanov
homologies $Kh(L)$ are invariant in the case of diagrams for which
the corresponding atoms are orienatble. The positive answer to
this question is given in the next section.

\section{Atoms and the Khovanov Complex \\ of 2-Fold Covering}

The main goal of the present section is to prove the following

\newcommand{\fF}{{\bf F}}

\begin{thm} ѕусть $\fF$ --- поле.
Let $L,L'$  be two equivalent virtual diagrams such that the
corresponding atoms are orientable. Then
$Kh_{\fF}(L)=Kh_{\fF}(L')$.\label{glth}
\end{thm}

The basic idea of the proof is the following. For each virtual
diagram $L$, consider the corresponding atom $V(L)$. Then, we
shall use the {\em orienting covering techniques}. Namely, if the
atom $V(L)$ is orientable, we consider two disjoint copies of the
atom $V(L)$, otherwise we consider the atom ${\tilde V(L)}$ which
is the orineting double cover over $V(L)$. It is defined as the
usual double cover over the correspondence 2-surface, herewith the
preimage of the frame is a graph to be considered as the frame
with cells cololred accordingly. The obtained atom is thus either
two-component or one-component depending whether the initial atom
is orientable or not.

Denote the obtained atom by $V^{2}(L)$, and the corresponding
virtual link by  $K(V^{2}(L))$.

This construction can be treated as follows: we consider two
disjoint sets of of crossings of the atom with a given
$A$-structure, and then connect some crossings by edges.

Thus, for each virtual knot, we have defined its ``covered
version'':

\begin{equation} L\to V(L)\to V^{2}(L)\to
Kh_{\fF}(K(V^{2}(L)))\end{equation}

Starting from knot diagrams, this construction is described as
follows. Given a virtual link diagram $L$ with $n$ classical
crossings $X_{1},\dots,X_{n}$. These crossings are connected
somehow between each other. Thus we have a graph $\Gamma$ immersed
in ${\bf R}^{2}$. Each crossing $X_{i}$ has four tails
$X_{i1},X_{i2},X_{i3},X_{i4}$ numbered, say, clockwise; herewith
the tails are connected to other tails by branches corresponding
to edges of the atom. Suppose the edge $l_{j}$ connects some two
tails $X_{j_{1}j_{2}}$ and $X_{j_{3}j_{4}}$, whereas
$j_{2},j_{4}\in\{1,2,3,4\}$.

The diagram $K(V^{2}(L))$ is then constructed as follows. It
contains $2n$ crossings $X'_{1},\dots, X'_{n}, X''_{1},\dots,
X''_{n}$ to be connected by branches. Each branch $l_{j}$ of the
initial diagram has two images: $l_{j}^{1}$ and $l_{j}^{2}$. Each
of the two edges $l_{j}^{i}$ connects one tail $X'_{j_{1}j_{2}}$
or $X''_{j_{1}j_{2}}$ with another tail $X'_{j_{3}j_{4}}$ or
$X''_{j_{3}j_{4}}$. For each edge $l^{1}_{j}$, we have to
indicate, which types of tails it should connect ($X'$ or $X''$).
Here we have some description arbitrarity. The point is that
before describing the edges, we have got no natural vertex
ordering, that is we did not say which of the vertices $X'_{i}$
and $X''_{i}$ is the first, and which one is the second one. To
overcome it, let us consider a maximal tree $\Delta$ in the graph
$\Gamma$ and decree that all edges $l_{j}^{1}$ corresponding to
edges of $\Delta$ connect the tails $X'_{j_{1}j_{2}}$ and
$X'_{j_{3}j_{4}}$ (thence, the edges $l_{j}^{2}$ connect some
tails of type $X''_{j_{1}j_{2}}$ и $X''_{j_{3}j_{4}}$). Another
choice of the maximal tree would lead to a notation change: for
some pairs, $X'_{j}$ and $X''_{j}$ will permute. After this, the
connection rule for the remaining tails by edges $l_{i}^{1}$ and
$l_{i}^{2}$ goes as follows. Again, we have some freedom: showing
which pairs of tails are to be connected by an edge
$l_{i}^{\cdot}$ we shall not indicate which of the two edges
$l_{i}^{1}$ or $l_{i}^{2}$ is used here; the ``antipodal'' pair of
tails obtained by permutation $X'\longleftrightarrow X''$ would be
also connected by an edge; it is not important for the diagram
what the name of the edge is. Also, we shall not pay attention to
the disposition of edges $l_{i}^{\alpha}$ on the plane since the
class of virtual link does not depend on that.

Thus, we have fixed a spanning tree $\Delta\subset \Gamma$. Each
edge  $l_{j}$ not belonging to this edge generates a minimal cycle
on the subgraph $\Delta\cup l_{j}\subset \Gamma$. In the case when
this cycle is {\em good} (see below), we use the edge $l_{j}^{1}$
to connect the tails $X'_{j_{1}j_{2}}$ and $X'_{j_{3}j_{4}}$ and
the edge  $l_{j}^{2}$ will be used for connecting
 $X''_{j_{1}j_{2}}$ and $X''_{j_{3}j_{4}}$. In the case of bad
 cycle we connect by $l_{j}^{1}$ the tails $X'_{j_{1}j_{2}}$ and
$X''_{j_{3}j_{4}}$, whence the tails $X''_{j_{1}j_{2}}$ and
$X'_{j_{3}j_{4}}$ are connected by $l_{j}^{2}$. The notions of
{\em good and bad edges} result from orienting and non-orienting
cycles on the atom. The covering maps any orienting cycle to a
cycle and any non-orienting cycle to a path with different
endpoints $X'_{k},X''_{k}$. Now, we define the notion of good edge
(for edges not belonging to $\Delta$) together with the notion of
orienting cycle. {\bf A cycle is good if the number of its
transverse passings through crossings is even.}

It is easy to check that this definition of a good cycle agrees
with the definition coming from the atom.

Thus, we have defined the notion of orienting cycle and good edge
(for edges not belonging to the tree $\Delta$). Thus, we have
completely constructed the virtual diagram $K(V^{2}(L))$. Note
that this definition does not depend (modulo detour moves) on the
choice of the tree $\Delta$.

Moreover, by the atom $V^{2}(L)$, the ``double-cover'' link is
recovered modulo virtualisation (however, this does not change
Khovanov homologies); we have given an explicit way for
constructing the diagram $K(V^{2}(L))$ by the diagram $L$; this
way corresponds to some immersion of the frame of $V^{2}(L)$
(which preserves the $A$-structure)

It is easy to see that if we apply the detour move to the initial
diagram $L$, the diagram $K(V^{2}(L))$ would also be operated on
by a detour move. Besides, we have the following
\begin{lm}
If we apply a Reidemeister move to a diagram $L$, then the diagram
$K(V^{2}(L))$ will be operated on by the same move in two places.
\end{lm}

\begin{proof}
We shall denote links before and after the Reidemeister move by
$L,L'$, denote the frames of the corresponding atoms by
$\Gamma,\Gamma'$ and the corresponding ``covered diagrams'' by
${\tilde L}, {\tilde L}'$, respectively.

Each Reidemeister move is a transformation of a knot diagram
inside a certain domain; inside this domain, the diagram $L$ is
represented by a subdiagram $P$, whence $L'$ is represented by
some diagram $Q$. Herewith we have some tails $t_{1},\dots,t_{n}$
that connect the subdiagram $P$ (of $L$) or the subdiagram $Q$ (of
$L'$) with the remaining fixed part of the diagram. In the case of
the first Reidemeister move we have $n=2$, in the case of the
second Reidemeister move $n=4$, for the third Reidemeister move we
have $n=6$. For the diagrams ${\tilde L}$ and ${\tilde L}'$, each
tail  $t_{i}$ has a lift $t'_{i},t''_{i}$ to the covering diagram.
If the diagram $P$ (and $Q$) does not contain non-orienting cycle
then  we have two copies of $P$ on ${\tilde L}$ (and two copies of
$Q$ on ${\tilde L}'$). The claim of the theorem is that the
connections of these copies agree with tails. A coordinated and a
non-coordinated cases are shown in Fig. \ref{sgns}.

\begin{figure}
\centering\includegraphics[width=400pt]{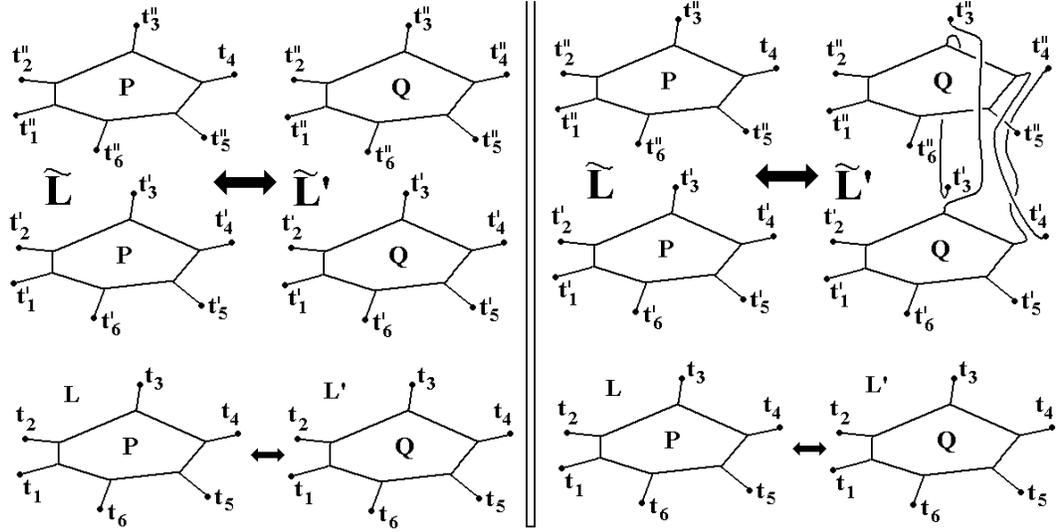}
\caption{Coordinated and non-coordinated connection} \label{sgns}
\end{figure}

In the coordinated case, all vertices $t_{i},t'_{i}$ are connected
coordinatedly for the lift of $L$ and for the lift of $L'$.
Namelt, all  $t'_{j}$ are connected to each other, and all
$t''_{k}$ are connected to each other. In the case shown in Fig.
\ref{sgns} the connection according to the lift of $L$ is as
before, whence the connection for the lift of $L'$ is distinct:
one copy of  $Q$ is connected to the tails
$t'_{1},t'_{2},t''_{3},t''_{4},t'_{5},t'_{6}$ whence the other one
connects the tails
$t''_{1},t''_{2},t'_{3},t'_{4},t''_{5},t''_{6}$. We have to show
that only the first type of connection takes place.

While applying the first Reidemeister move, we get a loop
consisting of one edge $l_{j}$ for which the initial point (the
same as the final point) divides some edge $l_{k}$. The edge
$l_{j}$ represents a good cycle. Besides, the edge $l_{k}$ is
divided into two edges  $l_{k_1}$ and $l_{k_2}$ which are not
opposite to each other in the separating vertex. Consider a cycle
$c$ of the graph $\Gamma'$ containing the edge $l_{k}$. With it,
one can naturally associate a cycle $c'$ containing $l_{k_1}$ and
$l_{k_2}$ on the graph $\Gamma'$. If  $c$ is orientable, then so
is  $c'$ and vice verse. Indeed, since the edges $l_{k_1}$ and
$l_{k_2}$ are not opposite, we have no additional transverse
point. Thus we have the same parity for the cycles $c$ and $c'$.
So, if the edge $l^{1}_{k}$ of the diagram $K(V^{2}(L))$ connects,
say, the tails  $X'_{pq}$ and $X''_{rs}$, then the edges
$l^{2}_{k1}l^{2}_{k2}$ connect the same tails  $X'_{pq}$ and
$X''_{rs}$. This means that the two pictures before and after the
Reidemeister move agree. Each of the two remaining Reidemeister
moves represents a reconstruction of some ``interior domain''
having some tails. For the second Reidemeister move, we have four
tails, for the third Reidemeister move, we have six tails.

Consider the second Reidemeister move. The diagram $L'$ contains a
bigon $cd$ and four emanating edges $a,b,e,f$, see Fig. \ref{2dR}.

\begin{figure}
\centering\includegraphics[width=200pt]{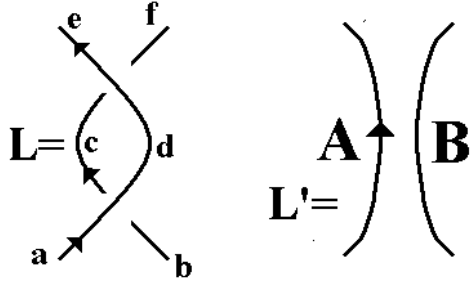} \caption{The
second Reidemeister move} \label{2dR}
\end{figure}

The bigon is evidently good. Thus, the whole set of edges can be
uniquely lifted to the diagram $K(V^{2}(L'))$. Finally we get two
sets of branches $a',b',c',d',e',f'$ and $a'',b'',c'',d'',e''$.
The branches $c',d'$ form a bigon which admits a decreasing
Reidemeister move. The same is true about the bigon $c'',d''$.
Applying the decreasing second Reidemeister moves to these bigons,
we get the diagram $K(V^{2}(L))$. Indeed, it suffices to show that
after such Reidemeister moves the edge $a'$ connects to the edge
$e'$, whence the edge $b'$ connects to $f'$ (herewith $a''$
connects to $e''$ and $b''$ connects to  $f''$). The last
statement follows from the fact that each cycle of  $L$ passing
through $a,c,e$ has the same number of transverse passings as the
corrsponding cycle  $L'$ passing through the edge $A$.

In the case of the third Reidemeister move, we have one hexagon
with six exterior branches for both $L$ and $L$'.

Both triangles represent good cycles, since they do not connect
transverse passes. Thus, the corresponding domains  $P$ and $Q$
are lifted to two copies of  $P$ and $Q$, respectively. We only
have to check that these two lifts agree.

To do it, we have to show that any two paths $\gamma\in L$ and
$\gamma'\in L'$ connecting  $t_{i},t_{j}$ have similar lifts to
$L,{\tilde L}$. For instance, if one preimage $\tilde \gamma$ of
the path $\gamma$ is lifted to the path connecting $t'_{i},t'_{j}$
and passing inside $P$ then each path $\gamma'$ from $Q$ having
the same points as $\gamma$ should have a preimage $\tilde\gamma'$
connecting $t'_{i},t'_{j}$ (and not $t'_{i},t''_{j}$). In Fig.
\ref{3rrd}, we give an example of two such paths.

\begin{figure}
\centering\includegraphics[width=240pt]{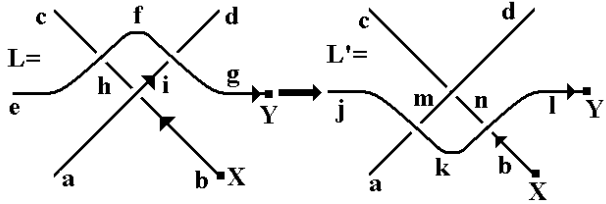} \caption{The
third Reidemeister move} \label{3rrd}
\end{figure}

Namely, consider some two paths between $X$ and $Y$ on the
diagrams $L$ and $L'$. We state, that their lifts  ${\tilde L}$
and ${\tilde L}'$ agree. This follows from the fact that the
number of transverse points for these paths is the same (equal to
zero). One can prove the analogous statement about paths, as shown
in Fig. \ref{3rrd}.
\end{proof}

Thus, by lemma \ref{lemm}, the Khovanov homologies of the
``doubled-covered'' knot are not changed while applying some
Reidemeister mvoes to the initial diagram. Thus we get the
following

\begin{thm}
The map $L\to Kh(K(V^{2}(L)))$ provides a well defined invariant
of virtual link.
\end{thm}

\begin{re} Note that a Reidemeister move may change the
orientability of an atom. Thus, for instance, for an orientable
atom $V(K)$ and the corresponding two-component atom $V^{2}(K)$, a
second Reidemeister move applied to $K$ may transform the atom
$V(K)$ to a non-orientable one, thus connecting the two components
of the atom $V^{2}(K)$.
\end{re}

Suppose the atom corresponding to a virtual link diagram $L$ is
orientable. Then $K(V^{2}(L))$ consists of two copies of $L$.
Since $\fF$ is a field, we have:
$Kh_{\fF}(K(V^{2}(L)))={Kh_{\fF}(L)}^{\otimes 2}$.

Thus, the homologies  $Kh(L)$ can be recovered from
$Kh(K(V^{2}(L)))$ by taking the tensor square root. In the case
when the coefficient ring is a field, we get a 2-variable
Poincar\'e polynomial with non-negative integer coefficients. So,
we have to extract the square root of this polynomial to get a
polynomial with non-negative coefficients. If it is possible, then
such a polynomial is unique. This completes the proof of Theorem
\ref{glth}. This leads to

\begin{thm}
Let  ${\bf F}$ be a field and let a link $L$ be such that
$Kh_{{\bf F}}(K(V^{2}(L)))$ can not be represented as a tensor
square. Then $L$ has no diagram for which the corresponding atom
is orientable. In particular, such a link $L$ can not be
classical.
\end{thm}

An important and interesting question is whether there exist two
non-isotopic classical links $L,L'$ that can be obtained from each
other by a sequence of generalised Reidemeister moves and
virtualisations. The Khovanov complex gives a partial answer to
this question.

From Theorem \ref{glth} and the invariance of Khovanov homologies
modulo virtualizations, we obtain the following result.

\begin{thm}
If a classical link  $L'$ can be obtained from another classical
link  $L$ by using generalised Reidemeister moves and
virtualisations, then $L$ and $L'$ have the same Khovanov
homologies with coefficients in any given field.
\end{thm}

{\bf e-mail: vassily@manturov.mccme.ru}

\end{document}